# Antimonopoly Regulation Method in Energy Markets Based on the Vickrey–Clarke–Groves Mechanism

Vadim Borokhov[1]

LLC "En+development", ul.Vasilisi Kozhinoi 1, Moscow, 121096, Russia

Abstract

**We evaluate the applicability of the generic Vickrey–Clarke–Groves (VCG) mechanism as an antimonopoly measure against a profit-maximizing producer with market power operating a portfolio of generating units at the centralized two-settlement energy market. The producer may indicate in its bid not only the altered cost function but also the distorted values of the technical parameters of its generating units, which enter the system-wide constraints of the centralized dispatch optimization problem. To ensure the applicability of the VCG method in this setting, we identify an additional assumption on the changes of the feasible set of the centralized dispatch optimization problem induced by variations of the producer's technical parameters. In the framework of the generic VCG mechanism, we propose an antimonopoly regulation method based on a regulator estimate of the producer's truthful bid. If this estimate is exact, the producer's maximum profit coincides with that in the case of the truthful bidding when no antimonopoly measure is applied. If the estimate is not exact, the error affects neither the producer's (weakly) dominant bid nor its optimal nodal output but manifests itself in the total uplift payment. This ensures an efficient allocation in the form of the optimal output/consumption schedule and shields the (pre-uplift) market prices from the producer's market power. We compare the suggested method with the alternative antimonopoly regulation approach based on the replacement of the producer's bid by a bid composed by the regulator.**

## I. INTRODUCTION

Many power sector reforms resulted in the formation of liberalized power markets with market-based pricing for power and security-constrained economic dispatch based on financially binding bids submitted by the market players (Bohn et al., 1984; Littlechild, 1988; Schweppe et al., 1988). The goal of the liberalized power market is to determine both the optimal output/consumption schedule of the market players and the market prices (and applicable side payments) that support this optimal dispatch. The latter is given by an optimal solution of the constrained optimization problem with an objective function usually having the form of the total social welfare of the market.

Exercise of market power in the liberalized power markets and estimate of its impact on the market have been the topics of intensive research (Joskow and Kahn, 2002; Wolak, 2003; Brennan,

[1] E-mail: vadimbor@yahoo.com. The views expressed in this paper are solely those of the author and not necessarily those of LLC "En+development".

2003; Wolak, 2005; Twomey et al., 2005; Brennan, 2006; Hesamzadeh et al., 2011), including the markets with forward contracts (Allaz, 1987; Allaz, 1992; Allaz and Vila, 1993; Hogan, 1997; Newbery, 1998; Green, 1999; Wolak, 2000; Bessembinder and Lemmon, 2002; Anderson and Xu, 2005). These works suggest that since the power systems have to account for both generator private constraints (such as minimum/maximum capacity limits, ramping constraints, minimal up/down times) and the power network constraints (the power flow equations, network power flow transmission constraints due to the thermal or security limits, etc.), the electricity markets proved to be prone to market power abuse even when no apparent dominant market players are present in the system and standard market power indicators have acceptable levels. Mount et al. (2001), Tierney et al. (2008), Fabra (2003) assessed different electricity market designs based on exposure to market power abuse by the market players.

A producer with market power may distort both the cost of power output and values of the technical parameters of its generating units, resulting in allocative inefficiency and economic surplus redistribution. A set of standard policies to reduce market power includes the measures that affect the market structure (such as forced divestiture of generating capacity, affiliation and merger control) and the behavioral limitations (introduction of mandatory hedging rate, application of tariffs or price-caps, etc.). Since the regulator lacks full information on the relevant cost components of a power producer (in particular, in the case of non-transparent fuel pricing) and the values of the technical parameters entering the producer private constraints on power output, an application of the behavioral limitations leads to the price and output/consumption volumes distortions relative to the ones obtained when the producer submits a truthful bid containing the truthful cost of power output and the truthful values of the technical parameters. Therefore, it is important to design an antimonopoly measure that minimizes the effect on the market prices and the output/consumption volumes from the error in the regulator’s estimate.

It is well known that a producer that can perfectly price discriminate all the consumers sells the same amount of goods as it does when it acts as a price-taker but (contrary to the price-taking supplier) captures the entire market surplus. Thus, it may prove to be beneficial for the market to create an economic environment, where a producer with market power is able to perfectly price discriminate the demand but is deprived of a share of the market surplus intended for the other market players. This approach is implemented in the Vickrey–Clarke–Groves (VCG) mechanism (Vickrey, 1961; Clarke, 1971; Groves, 1973), which is a generic truthful mechanism for achieving the allocative efficiency (which, for the case of power markets, means optimality of the resulting output/consumption schedule). The VCG mechanism aligns an agent’s profit with the total social welfare by using a special pricing method that produces a profit of each agent equal to the total social welfare (up to terms independent of the agent's strategy). Such a method was earlier developed by Clarke (1971), who proposed to calculate the agent's profit based on the difference in the total social welfare when a given agent participates in the market and when it does not. However, the VCG mechanism is not budget-balanced as the total cost paid by the buyers is generally lower than the total revenue received by the sellers, which limits the utilization of this approach in real auctions. Practical applicability of the VCG approach for real power markets (including the markets with non-convexities) has been extensively scrutinized by Hobbs et al. (1990), Rothkopf et al. (1990), Engelbrecht-Wiggans and Kahn (1991), Rothkopf and Harstad (1994), Rothkopf and Harstad (1995). These studies suggest that in practical applications the VCG mechanism may not be truth revealing due to the following. Since the real power markets are not isolated, the negative implications from the revelation of the private information to competitors may outweigh the benefits of truthful bidding at the power market. Also, the VCG mechanism is prone to collusive strategies that result in higher profits for the firms engaged in the exercise of market power. Bushnell and Oren (1994) researched the problems related to the truth revelation of the operation costs in electric power auctions. McGuire (1997) proposed to utilize the VCG algorithm for the unit commitment problem, while MacKie-Mason (1994) studied application of the VCG mechanism for power systems with transmission constraints. Moreover, the budget balancing problem of the VCG mechanism implies the associated uplift distributions, which may reduce the efficiency of the resulting outcome, (Krishna and Perry, 1998). A power market mitigation approach based on the generic VCG mechanism embedded in the form of the vesting contract was suggested by Borokhov (2010, 2011) for the case of a firm operating a single generating unit under the assumption that truthful values of the technical parameters of the firm's generating units are available to the regulator.

In this paper, we further research the applicability of the VCG mechanism in general two-settlement power markets as an antimonopoly measure against a profit-maximizing firm with market power operating a portfolio of the generating units located in possibly different nodes of the power system. The producer may manipulate not only the cost function of output but also the values of the technical parameters of its generating units. The latter may affect a set of possible market outcomes (i.e., the feasible set of the centralized dispatch optimization problem), which was usually assumed to be fixed in the previous studies. In particular, if these technical parameters enter not only the private constraints of the market player but also the system-wide constraints of the centralized dispatch optimization problem, the distortion of the technical parameters may have some nontrivial implications for the market efficiency. The application of the VCG approach only to market players subject to the antimonopoly regulation allows circumventing the budget-balancing problem at the cost of having the additional contribution to the total uplift introduced in the market. Based on the generic VCG algorithm, we propose an antimonopoly method that implies modifying the producer revenue function to ensure that its profit is a sum of the total social welfare and the terms that depend on the regulator estimate of the producer truthful bid and, therefore, are immune to the producer's strategy. If the regulator's estimate is not exact, the error may alter the total market uplift but affects neither the firm's optimal bid nor the optimal output volumes. Thus, the error has no influence on the planned dispatch schedule but manifests itself in the final prices through the total uplift payment for the market, which implies a redistribution of the total market surplus among the market players (with the total market surplus unaffected).

The paper is organized as follows. In section II we outline the VCG mechanism. The application of the VCG mechanism to a profit-optimizing producer at the two-settlement electric power market is formulated in Section III. Sections IV-V provide descriptions of the proposed antimonopoly algorithm at the day-ahead market (DAM). Section VI contains an example illustrating application of the proposed method to the firm with market power. The conclusions are presented in section VII.

Since most practical optimization problems in power markets involve continuous objective functions optimized over nonempty compact feasible sets, we assume that all maxima stated in the paper are attainable.

## II. THE VICKREY–CLARKE–GROVES MECHANISM

The VCG mechanism provides a method of revealing the individual preferences when they are not publicly observable and achieving efficient allocation of the resources. It is rather general and is applicable to the cases with non-convex bids. Consider a set of market players $I$ indexed by $i = 1,..,|I|$ participating in an auction where some commodity is traded. (The method is also applicable to a basket of commodities, but to simplify the notations we will consider just one good.) Each market player $i$ has the utility function $U_i(z_i, B_i)$, where $z_i$ denotes the amount of the commodity bought/sold by the market player and $B_i$ denotes the market player's bid containing information on the private constraints as well as the output cost (consumption benefit) function. If $i$ is a consumer (producer), then $U_i(z_i, B_i)$ corresponds to the (negative) value of the declared consumption benefit (output cost) function at $z_i$. For each $i$, the market player's truthful bid $B_i^{true}$ is assumed to be unknown to the other market players and the auctioneer. The auctioneer collects the bids $B_i$, $i = 1,..,|I|$, which potentially differ from the corresponding $B_i^{true}$, and determines the amount of commodity bought/sold by each market player as an optimal solution of the optimization problem

$$\max_{z \in S} \sum_{i \in I} U_i\,(z_i, B_i) \qquad (1)$$

with the feasible set $S$ describing the possible choices. The set $S$ is defined using both the private constraints of the market players and the system-wide constraints (such as the commodity balance constraint). Let $\{B\}$ designate the set of all bids, $\{B_{\bar{\iota}}\}$ stand for all the bids with $B_i$ excluded, $\bar{I}$ denote a set of all market players excluding $i$, and $z_1^*(\{B\}), .., z_{|I|}^*(\{B\})$ denote an optimal solution of (1). In this setting, the revenue of a supplier $i$ from delivering $z_i^*(\{B\})$ is given by $\sum_{j \in \bar{I}} U_j\big(z_j^*(\{B\}), B_j\big) + h_i(\{B_{\bar{\iota}}\})$ with some function $h_i(\{B_{\bar{\iota}}\})$ independent of the supplier bid. Similarly, for a buyer $i$ the expense function for consuming $z_i^*(\{B\})$ is given by $-h_i(\{B_{\bar{\iota}}\}) - \sum_{j \in \bar{I}} U_j\big(z_j^*(\{B\}), B_j\big)$. Following Mas-Colell et al. (1995), we show that in this setting the strategy $B_i^{true}$ of the market player $i$ is weakly dominant. If a supplier $i$ submits a bid $B_i$, then its profit at the optimal solution of (1) is given by

$$\pi_i(B_i) = \sum_{j \in \bar{I}} U_j\big(z_j^*(\{B_i, B_{\bar{\iota}},\}), B_j\big) + h_i(\{B_{\bar{\iota}}\}) + U_i(z_i^*(\{B_i, B_{\bar{\iota}},\}), B_i^{true}),$$

with the other market player bids $\{B_{\bar{\iota}}\}$ treated by the supplier $i$ as fixed (unknown) external parameters. Since $z_i^*(\{B_i^{true}, B_{\bar{\iota}}\})$ is an optimal point of the optimization problem $\max_{z \in S}\{\sum_{j \in \bar{I}} U_j\big(z_j, B_j\big) + U_i(z_i, B_i^{true})\}$, for any $z_i^*(\{B_i, B_{\bar{\iota}}\})$ we have

$$\sum_{j \in \bar{I}} U_j\big(z_j^*(\{B_i^{true}, B_{\bar{\iota}}\}), B_j\big) + U_i(z_i^*(\{B_i^{true}, B_{\bar{\iota}}\}), B_i^{true}) \geq \sum_{j \in \bar{I}} U_j\big(z_j^*(\{B_i, B_{\bar{\iota}}\}), B_j\big) + U_i(z_i^*(\{B_i, B_{\bar{\iota}}\}), B_i^{true}),$$

which implies that $\pi_i(B_i^{true}) \geq \pi_i(B_i)$. Thus, if a market player $i$ participates in the auction, then $B_i^{true}$ is its weakly dominant bid. If the participation in the auction is mandatory, then for any choice of the function $h_i(\{B_{\bar{\iota}}\})$, the market player's (weakly) dominant strategy is to submit its truthful bid. (We note that in this case a choice of $h_i(\{B_{\bar{\iota}}\})$ provides a price signal for the long-term planning, such as investment decisions of existing participants and potential new entrants.) However, if the participation in the auction is voluntary, then the function $h_i(\{B_{\bar{\iota}}\})$ has to satisfy certain additional participation constraints, such as a condition that the market player's profit is no lower than the highest profit it may receive should it withdraw from the auction. A special case of the VCG mechanism was earlier discovered by Clarke (1971) with $h_i(\{B_{\bar{\iota}}\})$ given by $h_i(\{B_{\bar{\iota}}\}) = -\max_{\substack{z \in S \\ z_i = 0}} \sum_{j \in \bar{I}} U_j\,(z_j, B_j)$, which is the optimal value of (1) when the market player $i$ is excluded from participating in the auction (it is assumed that the corresponding optimization problem is feasible). We note that the Clarke's choice for $h_i(\{B_{\bar{\iota}}\})$ generally implies higher profits for the firms with more market power, paying higher reward for the truthful bidding, which results in higher revenue (lower expense) per unit sold (bought) at the auction. Also, the profits of the market players are generally not invariant under the bid (dis)aggregation (i.e., if any two firms are considered as one unified firm, then the profit of the resulting firm may not be equal to the sum of the profits of the two original firms, and vice versa.)

The key shortcoming of the VCG mechanism is that in many cases (and under rather general assumptions) in the voluntary participation setting no choice of the functions $h_i(\{B_{\bar{i}}\}), i = 1, .., |I|$, ensures both the truthful bidding (i.e., the truthful bidding is a weakly dominant strategy for each market player) and the budget balancing (i.e, the total payment of the buyers equal the total revenue of the sellers), (Green and Laffont, 1979; Myerson and Satterthwaite, 1983).

We note the following specific features of the VCG mechanism. First, the objective function of the optimization problem (1) is additively separable into terms each depending on the bid of just one market player. Second, the feasible set $S$ (a.k.a. the choice set) is assumed to be unaffected by the difference between $B_i$ and $B_i^{true}$, $\forall i = 1, .., |I|$.

Applicability of the VCG mechanism to the power market was extensively scrutinized and criticized by a number of authors (Hobbs et al., 1990; Rothkopf et al., 1990; Engelbrecht-Wiggans and Kahn, 1991; Rothkopf and Harstad, 1994; Rothkopf and Harstad, 1995). In this paper, we emphasize the following additional issues. In a two-settlement power market, the planned output scheduled at the DAM can be amended at the real-time market resulting in the difference between the actual and planned outputs. If there is a statistical arbitrage between the DAM and the real-time market, then the supplier may bid not the truthful production cost but a blend with the real-time market expenses to redeem some of its DAM obligations. Furthermore, the technical parameters that are part of the supplier's bid may enter not only the producer private constraints but also the system-wide constraints. Theoretically, this may lead to a situation when a producer has economic incentives to distort the technical parameters of its generating units to affect a set defined by the system-wide constraints and, as a result, to manipulate the feasible set of (1). In these cases, the optimization problem (1) may not be equivalently formulated as the optimization problem with the same bids of the other market players but some modified bid of the producer in question with the truthful values of the producer's technical parameters (but a possibly different output cost function). In these instances, some additional assumption is needed to ensure applicability of the VCG approach.

## III. APPLICATION OF THE VCG MECHANISM TO POWER MARKETS

Consider a two-settlement electric power market (DAM and the real-time market) operating on bid-based security-constrained economic dispatch principle based on the financially binding bids submitted by the market players. DAM may have either power output/consumption volumes optimized separately from the other relevant products and services (such as ancillary services) or optimized jointly. Also, the unit commitment schedule can be either set before the DAM calculation or integrated into DAM optimization. Let DAM be based on the optimization of the total social welfare function $U$ with optimization variables $z$, which may include both discrete and continuous variables (as it is in the case of DAM simultaneously solving the unit commitment as well as the economic dispatch problem) or continuous variables only. A set of the explicit constraints, which includes both the market player private constraints and the system-wide constraints, is denoted as $\{C\}$. If DAM determines the unit commitment schedule as well, then the generating unit private constraints also include generator minimal up/down time limits, the relations between possible output volumes and the unit states, information on must run statuses, etc. We consider DAM with the hourly locational (nodal) marginal prices and the applicable side payments and uplifts. To avoid ambiguity, we clarify that the location marginal prices equal the Lagrange multipliers associated with the corresponding nodal power balance constraints and do not include any side payment/uplift distribution components. (Our analysis below is rather general and is applicable to any market pricing method that produces the power prices based on the power system parameters and the market player bids.)

The producer's DAM bid includes both the cost function and values of the technical parameters. We also suppose that the structure of power bids for DAM allows the producers to specify exactly all the relevant cost components of power output (such as power output costs, start-up, and no-load costs) and values of the technical parameters entering the generating unit private constraints. Let us denote as $\xi$ the set of technical parameters entering the constraints $\{C\}$ referring to both the supple side (the generating units' output limits, maximum ramp rates, etc.) and the demand side (for example, minimum/maximum consumption volumes) as well as the electric power network (such as power line admittances, maximum transmission capacities). The DAM optimization problem has the form

$$\max_{\substack{z\in M \\ s.t.\{C\}}} U(z) \qquad (2)$$

with the total social welfare function $U = U^{cons} + U^{prod}$, where $U^{cons}$ and $U^{prod}$ are the total benefit of power consumption as bid by the consumers and negative of the total cost of power output (and all the applicable products and services considered in DAM) as bid by the producers, respectively. If the consumers submit only perfectly inelastic DAM bids, then the term $U^{cons}$ is omitted from the function $U$. A feasible set of (2) denoted as $S(\xi)$ is assumed to be nonempty and compact. In what follows, we suppose continuity of the function $U$ in the non-discrete optimization variables.

We also assume that DAM and the real-time market converge and no arbitrage of any kind is possible for a producer between DAM and the real-time market: it is not profitable for the firm to sacrifice (a part of) its DAM earnings for (expected) complementary gain in the real-time market, and vice versa. Consequently, the profit maximization problem for a power producer at the wholesale power market cascades into the sequential solutions of the corresponding problems for DAM and the real-time market.

Consider a profit-maximizing firm $G$ operating a portfolio of generating units (possibly assigned to the different nodes of the power system) and selling the power output directly at DAM at the locational marginal prices (and the applicable side payments), i.e., not having any physical/financial contracts for power. Let variables $x$ be a subset of $z$ referring to the firm $G$ generating units, and variables $y$ denote the rest of variables $z$, so that $z = (x, y)$. Let $M_x$ and $M_y$ denote the corresponding product of the relevant finite discrete space (if any) and Euclidean space for variables $x$ and $y$, respectively. For example, in the case of DAM with integrated unit commitment problem the variables $x$ include both the output volumes and the binary state variables of the firm's generating units. Likewise, we denote as $\xi_G$ the technical parameters of the generating units of $G$ as specified by the firm and as $\xi_{\bar{G}}$ the rest of the technical parameters, so that $\xi = (\xi_G, \xi_{\bar{G}})$. Let $\{C_G(x, \xi_G)\}$ denote a set of private constraints of the firm's generating units, then $\{C\} = (\{C_G(x, \xi_G)\}, \{C_{\bar{G}}(z, \xi)\})$ with $\{C_{\bar{G}}(z, \xi)\}$ referring to the rest of $\{C\}$ and possibly depending on both $\xi_{\bar{G}}$ and $\xi_G$ (as it is in the case, for example, of the reserve adequacy requirements included in the DAM optimization). Let $S_x(\xi_G)$ be the private feasible set of $G$, i.e., a set of all $x \in M_x$ satisfying $\{C_G(x, \xi_G)\}$. For a given $x$ (not necessarily feasible), let us introduce a set $D_y(x, \xi)$ defined as a set of all $y$ satisfying $y \in M_y$ and the constraints $\{C_{\bar{G}}(x, y, \xi_G, \xi_{\bar{G}})\}$. The objective function of (2) is given by

$$U(x, y) = U_{\bar{G}}(y) - O_G(x), \quad (3)$$

where $O_G(x)$ denotes the sum of the declared cost functions of the firm $G$ generating units and $U_{\bar{G}}(y)$ describes the other market players stated benefit/costs functions. The function $O_G(x)$ is defined on $S_x(\xi_G)$. (To simplify the notations, we refer to $O_G(x)$ as the output cost function even if the producer sells at the DAM some other products and services in addition to electric power. Our consideration is rather general and includes these cases as well.) We note that the function $U(x, y)$ and, hence, $U_{\bar{G}}(y)$ may be independent of some of the variables $y$ as it is for the case of AC power flow model with variables $y$ including the voltage magnitude and the phase angle variables, which usually do not explicitly enter the total social welfare function. Also, assumed independence $U_{\bar{G}}(y)$ of the external parameters $\xi_{\bar{G}}$ doesn't imply a loss of generality as the set $y$ can be always enlarged to include additional optimization variables with the values constrained equal to $\xi_{\bar{G}}$. Analogously, introduction of $\xi_G$-dependence of $O_G(x)$ can be straightforwardly accounted for in the analysis below.

Let $B_G$ denote the firm $G$ DAM bid including both $O_G(x)$ and $\xi_G$. The truthful power output cost and the truthful values of technical parameters of the firm's generating units are denoted by $O_G^{true}(x)$ and $\xi_G^{true}$, respectively, with $O_G^{true}(x)$ defined on $S_x(\xi_G^{true})$. Hence, the firm's truthful bid is given by $B_G^{true} = \{O_G^{true}(x)\,, \xi_G^{true}\}$. If $B_G$ includes $\xi_G$ such that $S_x(\xi_G) \not\subset S_x(\xi_G^{true})$, then the system rebalancing might be needed in the real-time market because the firm $G$ is not able to deliver $x \notin S_x(\xi_G^{true})$. Therefore, to assess the market impact of such a solution, the real-time market analysis is needed. To circumvent this issue, we recall our assumption of no arbitrage between DAM and the real-time market, which for the case in question implies that for a solution with $x \notin S_x(\xi_G^{true})$ the expected real-time market penalty for the firm definitely outweighs its possible DAM extra gain. Therefore, we may restrict our consideration to a set $\Sigma_G = \{\xi_G | S_x(\xi_G) \subset S_x(\xi_G^{true})\}$ and bids $B_G$ with $\xi_G \in \Sigma_G$. We note that the function $O_G^{true}(x)$ is well-defined on $S_x(\xi_G)$, $\forall \xi_G \in \Sigma_G$. If $B_G$ differs from $B_G^{true}$, then the DAM outcome may result in the lower value of the total social welfare function $U^{true}(x, y) = U_{\bar{G}}(y) - O_G^{true}(x)$ and produce the corresponding deadweight loss for the market.

Let $z_0^* = (x_0^*, y_0^*)$ be an optimal point of (2) if the firm submits $B_G^{true}$ (if the DAM problem (2) has multiple maximizers, then $z_0^*$ denotes any of them). We have

$$U^{true}(z_0^*) = U_{\bar{G}}(y_0^*) - O_G^{true}(x_0^*), \text{ with } U_{\bar{G}}(y_0^*) = \max_{y \in D_y(x_0^*, \xi_G^{true}, \xi_{\bar{G}})} U_{\bar{G}}(y). \tag{4}$$

Since the economic goal of a profit-maximizing firm $G$ with market power can be misaligned with that of the market, i.e., optimization of the total social welfare $U^{true}(z)$, the bid $B_G$ may differ from $B_G^{true}$. As a remedy for that problem, we consider applying the VCG mechanism, which matches the profit maximization problem of $G$ with the DAM total social welfare maximization problem and forms economic incentives for the firm to submit the truthful bid. The stated no-arbitrage principle yields both that the firm $G$ attempts to maximize its DAM profit and that the firm's bid contains $\xi_G$ that satisfies $\xi_G \in \Sigma_G$. Let $z^* = (x^*, y^*)$ denote a solution to (2) when the firm $G$ submits a bid $B_G$. In the VCG mechanism, the revenue of $G$ is given by $\sum_{j\in\bar{G}} U_j\big(z_j^*(\{B\}), B_j\big) + h_i(\{B_{\bar{i}}\})$. In our framework, $\sum_{j\in\bar{G}} U_j\big(z_j^*(\{B\}), B_j\big) = U_{\bar{G}}(y^*(\{B_G, B_{\bar{G}}\}))$. In power market with additively separable total social welfare function, the only way the producer's bid $B_G$ may affect $y^*(\{B_G, B_{\bar{G}}\})$ is through $x$ and $\xi_G$, entering the system-wide constraints. Therefore, $y^* = y^*(\{B_{\bar{G}}\}, x^*(\{B_G, B_{\bar{G}}\}), \xi_G)$. As a result,

$$\sum_{j\in\bar{G}} U_j\big(z_j^*(\{B\}), B_j\big) = U_{\bar{G}}(y^*(\{B_{\bar{G}}\}, x^*(\{B_G, B_{\bar{G}}\}), \xi_G)) = \max_{y\in D_y(x^*(\{B_G, B_{\bar{G}}\}),\xi)} U_{\bar{G}}(y).$$

For any given $\xi_G$, varying the bid $B_G$ (having $\xi_G$ fixed), one can attain any feasible $x$ (which is defined as $(x,y)\in S(\xi)$ for some $y$ or, equivalently, $x\in S_x(\xi_G)$ provided that $D_y(x,\xi)\neq\emptyset$) as an optimal solution of (2). For example, this can be realized by choosing the cost function $O_G$ sufficiently low for the output below given $x$ and sufficiently high for the outputs above $x$. Consequently, instead of optimizing $B_G = \{O_G(x), \xi_G\}$ to achieve the maximum producer's profit one can use $x$ and $\xi_G$ as the optimization variables. Therefore, in the VCG mechanism the revenue of the producer $G$ as a function of $x$ and $\xi_G$ is given by

$$\max_{y\in D_y(x,\xi)} U_{\bar{G}}(y) + h_G(\{B_{\bar{G}}\}) \qquad (5)$$

with some function $h_G(\{B_{\bar{G}}\})$, which has no explicit or implicit dependence on the firm's bid (in particular, the latter implies that $h_G(\{B_{\bar{G}}\})$ has no dependence on $z^*$). The function $h_G(\{B_{\bar{G}}\})$ may depend on the bids $\{B_{\bar{G}}\}$, submitted by the other market players, and some other variables and parameters not influenced by $G$. For $\forall x\in S_x(\xi_G^{true})$, the DAM profit of the firm $G$ as a function of $x$ and $\xi_G$ is given by

$$\pi_G(x,\xi_G) = \max_{y\in D_y(x,\xi)} U_{\bar{G}}(y) + h_G(\{B_{\bar{G}}\}) - O_G^{true}(x), \qquad (6)$$

and the firm needs to find a profit-maximizing bid with some $\xi_G\in\Sigma_G$ and the cost function $O_G(x)$ defined on $S_x(\xi_G)$ so that the resulting $x^*$ maximizes its profit. (Formally, we put $\max_{y\in D_y(x,\xi)} U_{\bar{G}}(y) = -\infty$ if $D_y(x,\xi)=\emptyset$. However, this will not pose a problem since $D_y(x,\xi)=\emptyset$ entails that such $x$ cannot be a market outcome.) Since the firm is not able to influence the values of $\xi_{\bar{G}}$, we do not explicitly indicate $\xi_{\bar{G}}$ –dependence of $\pi_G$. We also note that the first term on RHS of (6) can be interpreted as the residual total social welfare for a given $x$. The profit-maximizing firm $G$ faces the problem

$$\max_{\substack{x,\xi_G \\ s.t. x\in S_x(\xi_G), \xi_G\in\Sigma_G}} \pi_G(x,\xi_G) = h_G(\{B_{\bar{G}}\}) + \max_{\xi_G\in\Sigma_G} \max_{(x,y)\in S(\xi_G,\xi_{\bar{G}})} U^{true}(x,y). \qquad (7)$$

To apply the VCG mechanism, we need to ensure that

$$\max_{\xi_G \in \Sigma_G} \max_{(x,y)\in S(\xi_G,\xi_{\bar{G}})} U^{true}(x,y) = \max_{(x,y)\in S(\xi_G^{true},\xi_{\bar{G}})} U^{true}(x,y). \quad (8)$$

In the general case, the feasible set of (2) depends on the firm $G$ generating units' technical parameters through both the private constraint set of $G$ and the system-wide constraints. In some cases, the intentional distortion of the technical parameter values by the firm can be equivalently represented as the modification of the firm's bid to include the truthful values of its technical parameters and the cost function inflated for $x$ such that $(x,y) \notin S(\xi_G,\xi_{\bar{G}})$, $\forall y,$ to signal that some points of $S(\xi_G^{true},\xi_{\bar{G}})$ are unattainable. For example, consider a uninode one-period power market with fixed demand. In this case, the technical parameters of the system include, among others, the generating unit maximum output limit. The constraint set of (2) is given by the private constraints of the market players and the power balance constraint. Hence, if the firm declares the unit's maximum output limit $g_{max}$ lower than its truthful maximum output limit $g_{max}^{true}$, then the producer's bid can be equivalently formulated using the truthful maximum output limit and the cost function properly inflated for the power output in the range $[g_{max}, g_{max}^{true}]$. This cost function modification ensures both that the objective function in (2) is still additively separable and that the certain points of $S(\xi_G^{true},\xi_{\bar{G}})$ are unattainable as optimal solutions to (2), which is equivalent to introduction of the additional constraints in the feasible set of (2). However, in some cases, the bid with distorted values of the technical parameters entering the optimization problem (2) cannot be equivalently formulated as a bid with the truthful values of these parameters and some modified cost functions retaining the additive separability of the objective function in (2). Therefore, we formulate the additional assumption that distortion of the firm's technical parameters cannot enlarge the feasible set of (2):

$$S(\xi_G,\xi_{\bar{G}}) \subset S(\xi_G^{true},\xi_{\bar{G}}) \text{ , } \forall \xi_G \in \Sigma_G. \quad (9)$$

If $\xi_G$ enters only the private constraint set of the firm (as it is in the case of the energy-only DAM with fixed unit commitment), then the condition (9) is satisfied. Indeed, in this case the constraint set $\{C_{\bar{G}}\}$, and, therefore, $D_y$ are independent of $\xi_G$. Using $S_x(\xi_G) \subset S_x(\xi_G^{true})$, $\forall \xi_G \in \Sigma_G$, we obtain $S(\xi_G,\xi_{\bar{G}}) = \cup_{x\in S_x(\xi_G)}[\{x\} \otimes D_y(x,\xi_{\bar{G}})] \subset \cup_{x\in S_x(\xi_G^{true})}[\{x\} \otimes D_y(x,\xi_{\bar{G}})] = S(\xi_G^{true},\xi_{\bar{G}})$, $\forall \xi_G \in \Sigma_G$. (We note that $\cup_{x\in S_x(\xi_G)}[\{x\} \otimes D_y(x,\xi_{\bar{G}})] = \{(x,y)|x \in M_x, y \in M_y, x \in S_x(\xi_G), y \in D_y(x,\xi_{\bar{G}})\}$.) Therefore, (9) holds. However, if $\xi_G$ enters the constraint set $\{C_{\bar{G}}(z,\xi)\}$ in some general way, the condition (9) may not be satisfied. As a theoretical example, consider a DAM with integrated unit-commitment procedure having a system-wide constraint (in addition to the power balance constraint) stating that the sum of the maximum output limits of the online units should not *exceed* a certain value. In such a system, a firm operating the generating unit with large maximum output limit $g_{max}^{true}$, low cost of output up to a certain level $g_L$, and high cost of output in the range $[g_L, g_{max}^{true}]$ may not be committed if it declares the truthful value of its maximum output limit. Meanwhile, the declaration of the intentionally distorted (lowered) value of the maximum output limit equal $g_L$ may enlarge the feasible set of (2) and result in the commitment of the unit in the optimal solution of (2). Clearly, (9) entails (8) due to the following.

$$\max_{\xi_G \in \Sigma_G} \max_{z \in S(\xi_G, \xi_{\bar{G}})} U^{true}(z) = \max_{z \in \cup_{\xi_G \in \Sigma_G} S(\xi_G, \xi_{\bar{G}})} U^{true}(z) = \max_{z \in S(\xi_G^{true}, \xi_{\bar{G}})} U^{true}(z).$$

Consequently, (7) yields

$$\max_{\substack{x, \xi_G \\ s.t.\ x \in S_x(\xi_G),\ \xi_G \in \Sigma_G}} \pi_G(x, \xi_G) = h_G(\{B_{\bar{G}}\}) + \max_{(x,y) \in S(\xi_G^{true}, \xi_{\bar{G}})} U^{true}(x, y) = \pi_G(x_0^*, \xi_G^{true}), \quad (10)$$

and the supplier's profit function attains its maximum value at $x = x_0^*$, $\xi_G = \xi_G^{true}$. We note that $\pi_G(x, \xi_G)$ may have the other maximizers as well. For example, if the set $D_y$ is independent of $\xi_G$, i.e., $D_y = D_y(x, \xi_{\bar{G}})$, and if $x_0^* \in S_x(\tilde{\xi}_G)$ for some $\tilde{\xi}_G \in \Sigma_G$, then $x = x_0^*$, $\xi_G = \tilde{\xi}_G$ is also the maximizer of $\pi_G(x, \xi_G)$. Similarly, a number of different bids may produce the same maximizer of $G$. Thus, the submission of the bid $B_G^{true}$ is a weakly dominant strategy. However, to compose the other optimal bid the firm $G$ may need some information on the bids of the other market players, which is unavailable to the firm. Therefore, the bid $B_G^{true}$ is the natural choice for the firm since no forecasting of unknown information is needed to compose it.

## IV. DESCRIPTION OF THE PROPOSED ANTIMONOPOLY REGULATION METHOD

Motivated by the VCG mechanism, we propose the following algorithm for a regulator to mitigate market power of the firm $G$.

a. After the gate closure time for the bid submission, DAM is calculated using the standard procedure. This results in DAM outcome $z^* = (x^*, y^*)$, the maximizer of (2), which is final for all the market players (including $G$). The resulting (pre-uplift) locational marginal prices are also final and are applied to all the market players excluding those subject to the antimonopoly regulation.
b. The regulator estimates both the power output costs and values of the technical parameters of the firm's generating units, which we denote as $O_G^{est}(x)$ and $\xi_G^{est}$ with $\xi_G^{est} \in \Sigma_G$, respectively. The firm's truthful bid estimated by the regulator is denoted as $B_G^{est} = \{O_G^{est}(x), \xi_G^{est}\}$. Since the regulator may not have full information needed to compose $B_G^{true}$, in the general case $B_G^{est}$ differs from $B_G^{true}$. The new DAM calculation is performed with the firm's DAM bid replaced by $B_G^{est}$ and using the original bids of the other market players.

Let us denote by $z_0^{est*} = (x_0^{est*}, y_0^{est*})$ a maximizer of (2) for this DAM calculation, we have $z_0^{est*} = z_0^{est*}(\{B_G^{est}, B_{\bar{G}}\})$. Also, let $R_G(x_0^{est*}, \xi_G^{est})$ be the revenue of $G$ from delivering $x_0^{est*}$ in DAM at the corresponding locational marginal prices (obtained in this new DAM calculation) with all the applicable DAM side payments and uplift payments associated with the DAM solution $z_0^{est*}$.

c. The regulator makes a decision on whether the firm should be subject to antimonopoly regulation in DAM. If the regulation is applied, the DAM revenue of $G$ is expressed as

$$R_G(x_0^{est*}, \xi_G^{est}) + U_{\bar{G}}(y^*) - U_{\bar{G}}(y_0^{est*}), \qquad (11)$$

which corresponds to (5) with $h_G(\{B_{\bar{G}}\}) = R_G(x_0^{est*}, \xi_G^{est}) - U_{\bar{G}}(y_0^{est*})$.

The stated algorithm should be known to the firm $G$ well in advance to be taken into account when composing its DAM bid. Thus, if $G$ is subject to the proposed antimonopoly regulation method, its DAM profit for $\forall x \in S_x(\xi_G^{true})$ is given by

$$\pi_G(x, \xi_G) = \max_{y \in D_y(x,\xi)} U_{\bar{G}}(y) - O_G^{true}(x) + R_G(x_0^{est*}, \xi_G^{est}) - U_{\bar{G}}(y_0^{est*}). \; (12)$$

If the participation in the centralized market is voluntary for the firm $G$, then it may find it beneficial to withdraw from the centralized market to avoid being subject to the pricing method (11) and sell its output through the bilateral contracts. Since a market player with market power has bargaining power, the contracted volumes and contract prices may differ from the ones obtained in the case of the truthful bidding. Therefore, we suggest mandatory participation of the producers with market power in the pricing algorithm “a” – “c”. In this case, the regulator needs to ensure that the revenue (11) covers the corresponding costs of the firms.

If the maximizer of (2) in the step “a” and/or the step “b” is not unique, then the regulator has to choose the maximizer in “b” that it considers appropriate. It is straightforward to verify that if the firm submits $B_G^{true}$, then the maximum value of the profit function $\pi_G(x, \xi_G)$ is independent of a choice for a maximizer in “a”. To the contrary, the maximum value of $\pi_G(x, \xi_G)$ generally depends on the choice of a maximizer in “b” if it is not unique. This reflects the fact that the various optimal points for (2) may result in the different profits for a given market player. Also, if the regulator estimate of $B_G^{true}$ is exact, i.e., $B_G^{est} = B_G^{true}$, and the firm submits its truthful bid, then both $z_0^*$ and $z_0^{est*}$ are maximizers of (2) with $U(z) = U^{true}(x, y)$, $\xi_G = \xi_G^{true}$, and we have

$$\max_{\substack{x, \xi_G \\ s.t.\; x \in S_x(\xi_G),\, \xi_G \in \Sigma_G}} \pi_G(x, \xi_G) = R_G(x_0^{est*}, \xi_G^{true}) - O_G^{true}(x_0^{est*}).$$

Thus, in the case of the exact regulator estimate of $B_G^{true}$, if (2) has multiple optimal points, then the maximum value of $G$'s profit is attained at the optimal point chosen by the regulator at the step "b". The proposed algorithm involves the following actions by the regulator: the estimate of the firm's truthful bid and the decision on whether to apply the antimonopoly measure. Since the complete information required for the exact estimate is unavailable to the regulator, the cases where the firm is underregulated or overregulated are inevitable. In these circumstances, the firm's DAM profit may differ from that obtained when $G$ submits its truthful bid to DAM and no antimonopoly actions are taken against the firm. For many other antimonopoly measures, this uncertainty results in distortions of both the volumes of goods delivered and the corresponding prices. However, if the firm behaves rationally, then the proposed method produces DAM schedule identical to that obtained when the firm submits its truthful bid $B_G^{true}$. Hence, the nodal output/consumption volumes in all nodes of the power system are shielded from market power of $G$. The natural question is whether the locational marginal prices are shielded as well. If the firm's offer is not marginal in any location of the firm's generating units, then the locational marginal prices in all nodes of the power system are protected from market power of the firm (provided that it submits a profit-maximizing bid). If the firm's bid is marginal in any of the nodes, then the locational marginal prices may be distorted by the firm's bid that is different from $B_G^{true}$ but still results in the profit maximizing nodal outputs. However, if the firm faces uncertainty in the values of its nodal output volumes resulting from DAM should it submit a different bid, then the firm's optimal bid is $B_G^{true}$. We also note that any firm's bid different from $B_G^{true}$ cannot yield additional profit but may lower it. Hence, the firm cannot profit from the intentional distortion of the locational marginal prices (unless it colludes with the other market players). Thus, if the firm's bid is not marginal or the firm faces uncertainty in the market outcome (provided that the firm does not collude with the other players), then all the locational marginal prices are also shielded from market power of $G$. Since uncertainty in the estimations of the firm $G$ truthful bid is objectively present, it results in the firm selling nodal power volumes at prices different from the corresponding locational marginal prices (and the applicable side payments) calculated at the step "a". This difference contributes to the total uplift, which can be allocated among the other market players in a way preserving the price signals and economic incentives produced by the market.

The first term on RHS of (12) coincides with the residual DAM social welfare calculated when the firm $G$ delivers $x$ and the other market players adjust their production/consumption volumes to maximize DAM total social welfare function treating $x$ as fixed. Thus, (12) (up to the terms that are independent of the firm's bid) coincides with DAM profit received by the firm in the case of perfect price discrimination of the residual demand in all hours of the day. The last two terms on RHS of (12) are independent of $x$ and do not alter the maximizer of $\pi_G(x, \xi_G)$ but are needed to deprive $G$ the share of the market surplus intended for the other market players (as forecasted by the regulator). The VCG mechanism has a clear economic meaning: it models the market environment that enables $G$ to perfectly price-discriminate the residual demand and provides economic incentives for the firm to submit the truthful bid to DAM. This aligns the economic interests of the firm with those of the market, i.e., the total social welfare maximization. Contrary to the Clarke's function $h_G(\{B_{\bar{G}}\})$, the proposed choice for $h_G(\{B_{\bar{G}}\})$ has no inherent rewarding of the market players with more market power for their truthful bidding but requires the regulator to estimate the market player's private information.

We also note that the antimonopoly regulation method formulated for a producer can be straightforwardly adapted for a consumer with market power, replacing the producer's revenue function (11) by the analogous expression for the consumer's expense function. Moreover, if the stated algorithm is to be applied to each firm in a group of market players $K$ (either producers or consumers or both) and the regulator has information that these firms have collusive strategies at the market, then $\forall i \in K$ the function $h_i$ should be independent of the bids of all the market players from the group $K$. This can be achieved, for example, by modifying the step "b" in the stated algorithm when applied to each market player $i$ from the group $K$ to have the bids of all the firms belonging $K$ (not just the bid of the market player $i$) replaced by the bids estimated by the regulator, resulting in the new DAM outcome $z_0^{est*} = z_0^{est*}(\{B_K^{est}, B_{\bar{K}}\})$, $\forall i \in K$, with $B_K^{est} = \{B_i^{est} | i \in K\}, \{B_{\bar{K}}\} = \{B_i | i \notin K\}$.

## V. EFFECT OF ERROR IN THE ESTIMATE OF THE FIRM TRUTHFUL BID

Exercise of market power by a firm results in distortion of the market surplus obtained by the firm and that of the other market players due to the reduction of the total social welfare and the surplus redistribution between the firm and the other market players given this (reduced) value of the total social welfare. Due to the ability of the firm to distort the market prices by varying its bid, the same decomposition takes place when the firm is subject to antimonopoly regulation based on the regulator's estimate of the firm's economic and technical aspects of power production because of the error in these assessments. For the comparison, we also consider the alternative antimonopoly regulation method that replaces the firm's bid by the bid composed by the regulator based on its estimate, $B_G^{est}$, - we will refer to such an algorithm as the "standard" regulation method. Thus, if the "standard" method is applied, then both of the abovementioned factors contribute, while if the proposed regulation method is used, then only the surplus redistribution occurs with no distortion of the optimal DAM schedule. Let us compare the profits received by $G$ in both methods. From (12) it follows that

$$\pi_G(x_0^*, \xi_G^{true}) = \pi_G^{st} + U^{true}(z_0^*) - U^{true}(z_0^{est*}), \qquad (13)$$

where $U^{true}(z_0^{est*}) = U_{\bar{G}}(y_0^{est*}) - O_G^{true}(x_0^{est*})$ is the total profit of the market players and $\pi_G^{st} = R_G(x_0^{est*}, \xi_G^{est}) - O_G^{true}(x_0^{est*})$ is the firm's profit when the "standard" regulation method is applied. From (9) and $\xi_G^{est} \in \Sigma_G$ we conclude that $U^{true}(z_0^*) \geq U^{true}(z_0^{est*})$, which implies $\pi_G(x_0^*, \xi_G^{true}) \geq \pi_G^{st}$. Thus, the firm's profit when subject to the suggested antimonopoly regulation method is at least as high as its profit in the "standard" regulation approach. Meanwhile, as it follows from (13), the total profit of the other market players is the same in both regulation methods:

$$U^{true}(z_0^*) - \pi_G(x_0^*, \xi_G^{true}) = U^{true}(z_0^{est*}) - \pi_G^{st}.$$

Thus, the proposed antimonopoly regulation method results in the same distortion of the total profit of the other market players as the "standard" regulation method and generates potentially higher profit for the firm subject to the regulation. This entails that the total profit of the market players (which equals the total social welfare) in the suggested regulation approach is no lower than that in the "standard" regulation method. Also, the proposed method produces incentives for the firm to exercise truthful bidding at the market. Hence, if the firm does not collude with the other players and behaves rationally, then the output/consumption volumes for all the market players and the locational marginal prices in all the nodes are unaltered by market power of the firm.

## VI. EXAMPLE

Let us illustrate the application of the suggested antimonopoly regulation method to a uninode one-period energy-only DAM with fixed unit commitment. Consider a power system with single producer $G$ operating one generating unit with the cost function $O_G^{true}(x) = \gamma x$, where $\gamma$ is some positive constant and $x$ is the unit's output satisfying the generating unit private constraints $\{C_G(x, \xi_G^{true})\} = \{0 \leq x \leq \xi_G^{true}\}$, where $\xi_G^{true}$ is the truthful value of the unit's maximum output limit. Let the aggregate consumer benefit function be given by $U_{\bar{G}}(y) = \alpha y - \beta y^2/2$ defined for the consumption volumes $y$, $0 \leq y \leq \xi_{\bar{G}}$, with $\xi_{\bar{G}} = \alpha/\beta$, where $\alpha$ and $\beta$ are some positive constants. The model parameters are assumed to satisfy $\frac{(\alpha-\gamma)}{2\beta} < \xi_G^{true} < \frac{(\alpha-\gamma)}{\beta}$, which also implies $\gamma < \alpha$ and $\xi_G^{true} < \xi_{\bar{G}}$. Thus, we have $\{C_{\bar{G}}(z, \xi_{\bar{G}})\} = \{0 \leq y \leq \xi_{\bar{G}}; x = y\}$. Since $\{C_{\bar{G}}(z, \xi_{\bar{G}})\}$ is independent of $\xi_G$, the assumption (9) is satisfied in this power system. If no antimonopoly regulation is applied, the firm's profit is given by $\pi_G^{m.p.}(x) = (\alpha - \beta x)x - \gamma x$. If the firm submits its truthful bid for DAM, then the maximum output limit constraint is binding. Hence, we have $x_0^* = y_0^* = \xi_G^{true}$ with

$$U^{true}(z_0^*) = (\alpha - \gamma)\xi_G^{true} - \beta(\xi_G^{true})^2/2.$$

In this case, the firm receives the profit $\pi_G^{m.p.}(x_0^*)$. If, however, the firm exercises its market power, then its profit-maximizing output is given by $x^{m.p.} = \frac{(\alpha-\gamma)}{2\beta}$. To achieve this market outcome, the firm may either distort its maximum output limit, i.e., indicate $\xi_G = x^{m.p.}$, and/or submit an inflated cost function to ensure that supply and demand intersect at $x^{m.p.}$, yielding the optimal solution $z^{m.p.} = (x^{m.p.}, y^{m.p.})$ with $y^{m.p.} = x^{m.p.}$. This results in the total social welfare $U^{true}(z^{m.p.}) = \frac{3(\alpha-\gamma)^2}{8\beta}$, which satisfies $U^{true}(z^{m.p.}) < U^{true}(z_0^*)$. Let the regulator estimate of the firm $G$ truthful bid be $B_G^{est} = \{O_G^{est}(x), \xi_G^{est}\}$ with $O_G^{est}(x) = \gamma^{est}x$, where $\xi_G^{est}$ and $\gamma^{est}$ are some positive constants assumed to satisfy $\xi_G^{est} < (\alpha - \gamma^{est})/\beta$, which entails $\gamma^{est} < \alpha$ and $\xi_G^{est} < \xi_{\bar{G}}$. In this case, $x_0^{est*} = y_0^{est*} = \xi_G^{est}$, and the corresponding marginal price equals $(\alpha - \beta\xi_G^{est})$. We also have $R_G(x_0^{est*}, \xi_G^{true}) = (\alpha - \beta\xi_G^{est})\xi_G^{est}$, and $h_G(\{B_{\bar{G}}\}) = -\beta(\xi_G^{est})^2/2$. We note that $D_y(x, \xi_{\bar{G}}) = \{y | y = x\}$ if $0 \le x \le \xi_{\bar{G}}$ and $D_y(x, \xi_{\bar{G}}) = \emptyset$, otherwise. Also, $S_x(\xi_G) = \{x | 0 \le x \le \xi_G\}$, $\Sigma_G = \{\xi_G | 0 \le \xi_G \le \xi_G^{true}\}$. Since $\xi_G^{true} < \xi_{\bar{G}}$, we have $\xi_G < \xi_{\bar{G}}$, $\forall \xi_G \in \Sigma_G$. Therefore, (12) implies $\pi_G(x) = \alpha x - \frac{\beta x^2}{2} - \gamma x - \frac{\beta\left(\xi_G^{est}\right)^2}{2} = \pi_G^{m.p.}(x) + \beta[x^2 - (\xi_G^{est})^2]/2$ for $x \in S_x(\xi_G)$, $\xi_G \in \Sigma_G$. Thus, if the firm is subject to the antimonopoly regulation, it faces the following profit optimization problem:

$$\max_{\substack{x, \xi_G \\ s.t.\ x \in S_x(\xi_G),\ \xi_G \in \Sigma_G}} \pi_G(x) = \max_{x \in [0, \xi_G^{true}]} \pi_G(x) = \pi_G^{m.p.}(x_0^*) + \beta[(\xi_G^{true})^2 - (\xi_G^{est})^2]/2. \quad (14)$$

It is straightforward to verify that the output volume $x_0^* = \xi_G^{true}$ maximizes (14). Therefore, if the proposed antimonopoly method is applied to the firm, then its output volume maximizing the firm's profit also maximizes the total social welfare function. In this case, the DAM output/consumption schedule and the marginal price, which equals $(\alpha - \beta\xi_G^{true})$, are unaffected by the producer's market power. The error in the regulator's estimate of the producer's bid manifests itself in the uplift (which can be of any sign) in the amount of $\beta[(\xi_G^{true})^2 - (\xi_G^{est})^2]/2$.

## VII. CONCLUSIONS

The electricity markets are prone to exercise of market power by producers, which reduces productive and allocative efficiency of the market, distorts price signals, and results in associated deadweight loss. Exercise of market power leads to distortions of both the nodal power output/consumption volumes and the locational marginal prices. The VCG mechanism is a generic truthful revelation mechanism for achieving an efficient allocation, which motivates applying this pricing method as an antimonopoly measure. However, the VCG mechanism is generally not budget balanced if applied to all market players in the setting of voluntary participation in the centralized market.

We studied the applicability of the VCG mechanism as an antimonopoly regulation method for a two-settlement power market and considered the case when a producer with market power may manipulate not only the cost function but also the technical parameters of its generating units. (In this setting, application of the VCG approach only to market players subject to the antimonopoly regulation allows circumventing the budget-balancing problem at the cost of having the additional contribution to the total uplift.) We showed that if these technical parameters enter the system-wide constraints, the additional assumption on the feasible set of the centralized dispatch optimization problem is needed to ensure the applicability of the VCG approach. In the framework of the generic VCG algorithm, we propose an antimonopoly measure that results in the producer's profit given by a sum of the total social welfare and the terms not influenced by the given producer. These terms depend on the bids submitted by the other market players and the regulator's estimate of the truthful bid of the producer. If the regulator's estimate is exact, the producer's profit equals its profit received when it submits the truthful bid and no antimonopoly regulation is applied to the firm. The method requires the regulator to estimate the producer power output cost function and the values of technical parameters of its generating units. Such estimates objectively involve a degree of uncertainty since the regulator may not have full information on the producer's economical and technical aspects of power production. The error in the regulator's estimate of the truthful bid has no influence on the centralized dispatch optimal solution (which follows from the allocative efficiency of the generic VCG mechanism) but contributes to the total uplift applied to the other market players. Therefore, the error has no effect on both the nodal output/consumption volumes and the locational marginal prices (which are the pre-uplift prices) as well as the firm's optimal bid but manifests itself in the final prices through the uplift payment allocation. This implies a redistribution of the total market surplus among the market players without affecting the total market surplus. Also, the proposed method has no inherent rewarding of the market players with more market power for their truthful bidding.

We compared the proposed antimonopoly regulation method with the alternative market power mitigation approach, in which the producer's bid is replaced by the bid composed by the regulator based on its estimate. The error in the estimate generally distorts both the resulting dispatch and the market prices. We showed that the total profit of the other market players is the same in both methods. Thus, in both methods, the error in the regulator estimate of the producer's truthful bid leads to the identical distortions of the total profit of the other market players. However, the proposed antimonopoly regulation algorithm, as opposed to the "standard" method, yields no reduction in the total social welfare and shields both the optimal dispatch and the (pre-uplift) locational marginal prices from the error in the regulator's estimate of the producer's truthful bid.

The suggested approach can be straightforwardly extended to the cases of a profit-maximizing consumer with market power and a group of the market players with collusive strategy at the market. We also note that the proposed antimonopoly method is independent of the specific pricing algorithm for power used in the market and is applicable for power markets with non-marginal pricing.